\documentclass[11pt]{article}
\usepackage{amsmath}
\numberwithin{equation}{section}
\usepackage{amsfonts}
\usepackage{amsthm}
\usepackage{amssymb, setspace}
\usepackage{color,graphicx,epsfig,geometry,hyperref,fancyhdr}
\usepackage[T1]{fontenc}
\usepackage{float}
\usepackage{subfig}
\usepackage{commath}
\usepackage{bbm}
\usepackage{mathrsfs,fleqn}
\usepackage{mathtools}
\newtheorem{theorem}{Theorem}[section]

\newcommand{\N}{\mathbb{N}}

\newcommand{\R}{\mathbb{R}}

\newcommand{\dnu}{\partial_\nu}

\newcommand{\grad}{\nabla}

\newcommand{\exD}{\mathbb{R}^2 \setminus \overline{D}}

\newcommand{\dd}{\,\mathrm{d}}
\newcommand{\pdr}{\partial_r}

\graphicspath{{paper-pics/}}

\begin{document}

\begin{flushleft}
\Large 
\noindent{\bf \Large Analytical properties of the principal clamped transmission eigenpair for a planar region}

\end{flushleft}

\vspace{0.2in}

{\bf  \large Haoran Qi and Isaac Harris\footnote{corresponding author \texttt{harri814@purdue.edu}} }\\
\indent {\small Department of Mathematics, Purdue University, West Lafayette, IN 47907, USA }\\
\indent {\small Email: \texttt{haoran.qi@purdue.edu} ;  \texttt{harri814@purdue.edu} }\\

{\bf  \large Andreas Kleefeld}\\
\indent {\small Forschungszentrum J\"ulich GmbH, J\"ulich Supercomputing Centre, Wilhelm-Johnen-Str., }\\
\indent {\small 52425 J\"ulich, Germany.}\\
\indent {\small FH Aachen - University of Applied Sciences, Faculty of Medical Engineering and }\\
\indent {\small Technomathematics, Heinrich-Mu\ss{}mann-Str. 1, 52428 J\"ulich, Germany.}\\
\indent {\small FH Aachen - University of Applied Sciences, Institute for Data-driven Technologies, }\\
\indent {\small Heinrich-Mu\ss{}mann-Str. 1, 52428 J\"ulich, Germany.}\\
\indent {\small Email: \texttt{a.kleefeld@fz-juelich.de}}\\


\begin{abstract}
\noindent This paper investigates the analytical properties of the principal eigenpair for the clamped transmission eigenvalue problem, which arises in the context of biharmonic (i.e. flexural wave) scattering by a clamped planar obstacle. By deriving a suitable Rayleigh quotient, we provide rigorous proofs for several fundamental properties that were previously only observed through numerical experiments. One of the results we have shown here is that the first eigenvalue is monotone decreasing with respect to the domain. Furthermore, we prove that the corresponding eigenfunction is of a single sign, a property analogous to the principal eigenfunction of the Dirichlet eigenvalues of the Laplacian. Another key result is a Faber--Krahn type inequality, establishing that the disk minimizes the principal eigenvalue among all domains of a given measure. The theoretical results are supported by a comprehensive numerical investigation.
\end{abstract}



\section{Introduction}\label{intro}
The inverse scattering problem for flexural waves (e.g. \cite{near-lsmBH,iterative-BH23,app2}), which are modeled by the biharmonic wave equation, has significant applications in the non--destructive testing of thin elastic plates. In particular, a central question that arises when applying {qualitative inversion methods} is the analysis of the so--called transmission eigenvalue problem. The eigenvalue problem we will study here comes from the scattering of an incident wave by an obstacle $D \subset \R^2$ that is a non--empty, simply connected bounded open set with $C^2$--boundary $\partial D$. In this paper, we are concerned with the specific case of a clamped planar obstacle, leading to the {\it clamped transmission eigenvalue problem}. While this problem has been investigated in previous studies \cite{HLK-TEP,HarrisKleefeldLee2025}, some of the fundamental analytical properties have only been studied numerically.

Here, the associated direct scattering problem is given as follows. For a given time--harmonic incident plane wave, denote $u^{\text{inc}}(x) = \text{e}^{\text{i}kx\cdot d}$, with incident direction $d \in \mathbb{S}^1$ (i.e., the unit circle) and wave number $k>0$, find the corresponding scattered field $u^{\text{scat}}$ that satisfies the 4th--order boundary value problem 
\begin{align}
\Delta^2 u^{\text{scat}} - k^4 u^{\text{scat}} = 0 \quad & \text{in } \mathbb{R}^2\setminus \overline{D}, \label{biharmonic} \\
u^{\text{scat}} \big|_{\partial D}=-u^{\text{inc}}  \quad \text{ and } \quad \dnu{u^{\text{scat}}} \big|_{\partial D} &=- \dnu{u^{\text{inc}}}. \label{cbc}
\end{align}
In the above boundary condition \eqref{cbc}, we let $\nu$ denote the outward unit normal vector on the boundary $\partial D$. To complete the system we impose that the scattered field satisfies the radiation condition  
\begin{align}\label{SRC}
\lim_{r \to \infty} \sqrt{r}(\pdr{u^{\text{scat}}} - \text{i} k u^{\text{scat}}) = 0  \quad  \quad r=|x|,
\end{align}
which is assumed to hold uniformly in any direction $\hat{x}=x/|x|$. It has been shown in \cite{biwellposed,DongLi24} that the direct scattering problem \eqref{biharmonic}--\eqref{SRC} is well--posed. 

Due to the operator splitting of the Biharmonic Helmholtz operator $\Delta^2  - k^4$ one can show that (see, e.g., \cite{DongLi24,HarrisLiOzochiawaeze2026}) has the following asymptotic behavior 
\begin{align*}
u^{\text{scat}}(x, d) = \frac{\text{e}^{\text{i} \pi/4}}{\sqrt{8\pi k}}\cdot \frac{\text{e}^{{\mathrm{i}k}|x|}}{\sqrt{|x|}}u^\infty(\hat{x},d) + \mathcal{O} \left( \frac{1}{|x|^{3/2}} \right), \quad \text{as } |x| \to \infty.
\end{align*}
Here, $u^\infty(\hat{x}, d)$ is the far field pattern or scattering amplitude of the scattered field. With this, one can define the associated far field operator 
$$F: L^2(\mathbb{S}^1) \to L^2(\mathbb{S}^1) \quad \text{ given by } \quad (Fg)(x) = \int_{\mathbb{S}^1} u^\infty (\hat{x}, d) g(d) \, \dd s(d).$$
The far field operator is the data that one uses when studying qualitative inversion methods like the linear sampling and factorization methods, see manuscripts \cite{Cakoni-Colton-book,kirschbook}. These methods can recover the obstacle $D$ with little a priori information. It is well--known that for these methods to be valid, we must have that the far field operator is injective with dense range. In \cite{GuoLongWu2024,HarrisLiOzochiawaeze2026,Zhu2026} it has been proven that $F$ is injective with dense range if and only if $k>0$ is not an associated {clamped transmission eigenvalue}. Therefore, this eigenvalue problem is instrumental in the shape reconstruction problem. 

The primary goal of this work is to provide a rigorous mathematical analysis of the principal clamped transmission eigenpair. For the clamped transmission eigenvalue problem, we get a homogeneous boundary value problem defined in all of $\R^2$. The existence, discreteness and determination via the far field operator was studied in \cite{HLK-TEP,HarrisKleefeldLee2025}. Here, we are concerned with proving analytical results for the associated eigenpair. The techniques used here are centered on the derivation of a suitable Rayleigh quotient, which allows us to study the properties of the eigenvalues and eigenfunctions through a variational framework. This approach will allow us to prove several key properties of the principal eigenpair. The main contributions of this paper are as follows. We establish that the principal clamped transmission eigenvalue is monotone decreasing with respect to the region $D$. This confirms that as the size of the scattering obstacle increases, the corresponding first eigenvalue decreases. We also show that the associated principal eigenfunction is of a single sign, a property analogous to the well--known behavior of the principal eigenfunction for the Dirichlet eigenvalue problem for $- \Delta$ in $D$. Finally, we derive a Faber--Krahn type inequality, which shows that the disk is the minimizer for the principal eigenvalue among all domains of a given measure. These theoretical results are further supported by numerical investigations.

The remainder of this paper is organized as follows. In Section \ref{RQ-analysis}, we formulate the biharmonic direct scattering problem and formally define the corresponding clamped transmission eigenvalue problem. In Section \ref{eigpair-analysis}, we introduce the equivalent variational formulation, construct the Rayleigh quotient, and prove our main analytical theorems regarding monotonicity, the Faber--Krahn inequality, and the sign property of the principal eigenfunction. In Section \ref{numerics}, we present numerical results substantiating our theoretical findings. Finally, we draw conclusions and discuss future research directions in Section \ref{fin}.

\section{Rayleigh Quotient for the Transmission Eigenvalues}\label{RQ-analysis}
In this section, we will derive an appropriate Rayleigh quotient for the transmission eigenvalue problem we wish to study that is associated with the aforementioned biharmonic scattering problem. Therefore, motivated by the previous discussion we define the clamped transmission eigenvalues as the values $k>0$ such that there is a non--trivial solution to the boundary value problem
\begin{align}
\Delta v + k^2 v = 0 \quad \text{in } D \quad &\text{and}  \quad   \Delta w - k^2 w = 0 \quad \text{in } \exD,  \label{tep1} \\
v+w = 0 \quad & \text{and } \quad  \dnu(v+w) = 0 \quad \text{on } \partial D \label{tep2}
\end{align}
with $(v, w)  \in H^1(D) \times H^1(\exD)$. See, e.g., \cite{HLK-TEP} for a derivation of the transmission eigenvalue problem from the associated direct scattering problem \eqref{biharmonic}--\eqref{SRC}. In our study, we let  $\nu$ denote the unit outward normal along $\partial D$. Here, we wish to study some of the analytical properties that were numerically investigated in \cite{HLK-TEP,HarrisKleefeldLee2025}. These are the dependence of the principal eigenvalue with respect to the Lebesgue measure of $D:=|D|$ as well as other analytical properties of the principal eigenpair. In particular, for other eigenvalue problems (see, e.g., \cite{robin-dirichletEig,steklov-robinEig,MagneticLaplacian-EigInequ,p-laplace} for recent studies) an appropriate Rayleigh quotient is often used to great success when studying eigenvalue problems associated with differential operators. 

Due to the structure of this eigenvalue problem, we will now analysis an appropriate Rayleigh quotient derived from \eqref{tep1}--\eqref{tep2}. To this end, notice that in \cite{HLK-TEP,HarrisKleefeldLee2025} it was proven that there is an infinite discrete set of positive clamped transmission eigenvalues. With this, we see that if $(v, w)$ satisfy \eqref{tep1}--\eqref{tep2} with $k>0$ then the real and imaginary parts also satisfy \eqref{tep1}--\eqref{tep2} with the same eigenvalue $k$. So we will assume that the corresponding eigenfunctions are real--valued. To aid in our study, we define an equivalent eigenfunction $u \in H^1(\R^2)$ as
\begin{equation}\label{new-eigfunc}
u = \begin{cases}
         \,\,\,\,\,\, v, \, &D, \\
        -w, \, &\exD.
    \end{cases}
\end{equation}
Therefore, the clamped transmission eigenvalue problem can be rewritten as 
\begin{equation} \label{eqn_trans_eig_transformed}
    -\Delta u = k^2(2\cdot\mathbf{1}_D - 1)u \quad \text{in } \R^2 
\end{equation}
where $\mathbf{1}_D$ is the indicator function on the set $D$. 
From the construction we see that the eigenvalue problems are equivalent, which is given in the following result.

\begin{theorem}\label{thm_eigequiv}
There exists a non--trivial eigenpair $k>0$ and $u \in H^1(\R^2)$ satisfying \eqref{eqn_trans_eig_transformed} if and only if there is a non--trivial eigenpair $k>0$ and $(v, w)  \in H^1(D) \times H^1(\exD)$ satisfying \eqref{tep1}--\eqref{tep2}.
\end{theorem}

With this in mind, we can use the equivalent variational formulation of \eqref{eqn_trans_eig_transformed} to study the clamped transmission eigenvalue problem. Therefore, we see that the equivalent variational formulation of \eqref{eqn_trans_eig_transformed} is given by 
\begin{equation} \label{varform-4RQ}
\int_{\R^2} \nabla u\cdot \nabla\phi \, \dd x = k^2\int_{\R^2} (2\cdot\mathbf{1}_D-1)u\phi \, \dd x \quad \text{for all} \,\,  \phi \in H^1(\R^2). 
\end{equation}
Note that one of the main difficulties is that this is an eigenvalue problem in all of $\R^2$ and the `potential' function $(2\cdot\mathbf{1}_D-1)$ is not sign definite. Now, since we have the existence of $k>0$ that are clamped transmission eigenvalues we let $u$ be the corresponding eigenfunction and take $\phi=u$ in \eqref{varform-4RQ} to obtain 
$$\int_{\R^2} |\nabla u|^2 \, \dd x = k^2 \left(2\int_D u^2\, \dd x - \int_{\R^2} u^2\, \dd x \right).$$
This implies that for any eigenfunction 
$$2 \| u \|^2_{L^2(D)} - \| u \|^2_{L^2(\R^2)} >0 \implies  2\| u \|^2_{L^2(D)} > \| u \|^2_{L^2(\R^2)} .$$ 
So to study this problem, similar to the analysis in \cite{mainref} we will assume that 
$$u \in H^1(\R^2)  \quad \text{such that } \quad 2\| u \|^2_{L^2(D)} > \| u \|^2_{L^2(\R^2)}$$
which is a non--trivial subset of $H^1(\R^2)$ since the extension by zero of any $C^\infty_c (D)$ function satisfies the above inequality. 

By the equivalent variational formulation, we can now define the Rayleigh quotient associated with \eqref{eqn_trans_eig_transformed} and equivalently \eqref{tep1}--\eqref{tep2} as 
\begin{equation}
    \label{eqn_frac_form}
    \mathcal{R}_D  [u] := \frac{\| \grad u \|^2_{L^2(\R^2)}}{ 2\| u \|^2_{L^2(D)}- \| u \|^2_{L^2(\R^2)} } \quad \text{such that } \quad 2\| u \|^2_{L^2(D)}- \| u \|^2_{L^2(\R^2)}>0,
\end{equation}
for $u\in H^1(\R^2)$. This quantity will be used to study the analytical properties of the principal eigenpair corresponding to the first (i.e. smallest) eigenvalue $k_1>0$. To this end, we wish to connect the minimizer of the Rayleigh quotient to the first clamped transmission eigenvalue. In the following result we show that the minimal value of $\mathcal{R}_D [u]$ exists and is equal to the first clamped transmission eigenvalue. 

\begin{theorem}\label{thm_first_eig_in_var_form}
Let $k_1$ denote the first clamped transmission eigenvalue of \eqref{eqn_trans_eig_transformed} or equivalently \eqref{tep1}--\eqref{tep2}, corresponding to the region $D$. Then we have that 
    \begin{equation}
        \label{eqn_inf_quotient}
        k_1^2 \,\, = \min\limits_{\substack{u\in H^1(\R^2)} } \left\{  \mathcal{R}_D  [u]  \,\, : \,\, 2\| u \|^2_{L^2(D)}- \| u \|^2_{L^2(\R^2)}>0 \right\}.
    \end{equation}
    Moreover, if $u_1$ is the corresponding eigenfunction, then it is the minimizer satisfying 
    \[
        k_1^2 = \mathcal{R}_D  [u_1]\leq \mathcal{R}_D  [u] \quad \text{for all $u \in H^1(\R^2)$ : } \,\, 2\| u \|^2_{L^2(D)}- \| u \|^2_{L^2(\R^2)}>0.
    \]
\end{theorem}
\begin{proof}
We begin by proving that the Rayleigh quotient is bounded below by a positive constant. To this end, note that by the constraint we have that 
$$\| u \|^2_{L^2(D)} \leq  \| u \|^2_{L^2(\R^2)}  < 2 \| u \|^2_{L^2(D)}.$$
Moreover, by the Cauchy--Schwarz inequality we obtain the estimate 
$$ \| u \|^2_{L^2(D)}  \leq  |D|^{1/2}\|u^2\|_{L^2(\R^2)} = |D|^{1/2}\|u\|_{L^4(\R^2)}^2.$$
By appealing to the Gagliardo--Nirenberg (Ladyzhenskaya) inequality given by 
$$ \|u\|_{L^4(\R^2)} \leq C\|\nabla u\|_{L^2(\R^2)}^{1/2}\|u\|_{L^2(\R^2)}^{1/2} \quad \text{for all $u \in H^1(\R^2)$ }$$
(see, e.g., Chapter 12 in \cite{SobSpaceBook}) we can further estimate to get that 
\begin{align}
        \label{funcspaceupbound}
        \| u \|_{L^2(D)} < C|D|^{1/2} \|\nabla u\|_{L^2(\R^2)} \quad \text{for all $u \in H^1(\R^2)$ : } \,\, 2\| u \|^2_{L^2(D)}- \| u \|^2_{L^2(\R^2)}>0.
\end{align}
Therefore, we can obtain the lower bound
$$ \mathcal{R}_D  [u] = \frac{\| \grad u \|^2_{L^2(\R^2)}}{2\| u \|^2_{L^2(D)}- \| u \|^2_{L^2(\R^2)}} \geq \frac{\| \grad u \|^2_{L^2(\R^2)}}{2\| u \|^2_{L^2(D)}} > \frac{C}{ |D|}.$$
This implies that the infimum of the Rayleigh quotient exists and note that the constant $C>0$ is independent of $D$. 

Now that we know that the infimum exists in the subset of $H^1(\R^2)$ satisfying the constraint, we consider the minimizing sequence
$$\{u_j\}_{j \in \N} \subset \big\{u \in H^1(\R^2)\,\,:\,\,  2\| u \|^2_{L^2(D)}- \| u \|^2_{L^2(\R^2)}>0 \big\}$$ 
such that 
$$\lim\limits_{j \to\infty}\mathcal{R}_D [u_{j}] = \inf\limits_{\substack{u \in H^1(\R^2)}} \left\{  \mathcal{R}_D [u] \,\, : \,\,   2\| u \|^2_{L^2(D)}- \| u \|^2_{L^2(\R^2)}>0 \right\}.$$
Therefore, by the constraint, we can `normalize' the minimizing sequence with the denominator of the Rayleigh quotient 
$$\tilde{u}_j =  \frac{u_j}{ \sqrt{2\| u_j \|^2_{L^2(D)}- \| u_j \|^2_{L^2(\R^2)}}  } \quad \text{ for all $j \in \N$}$$
which gives that 
$$\{\tilde{u}_j\}_{j \in \N} \subset \big\{u \in H^1(\R^2)\,\,:\,\,  2\| u \|^2_{L^2(D)}- \| u \|^2_{L^2(\R^2)}=1\big\}.$$ 
This would imply that 
$$\lim\limits_{j \to\infty} {\|\nabla \tilde{u}_j\|_{L^2(\R^2)}^2} \,\, = \inf\limits_{\substack{u \in H^1(\R^2)}} \big\{  \|\nabla {u}\|_{L^2(\R^2)}^2 \,\, : \,\, 2\| u \|^2_{L^2(D)}- \| u \|^2_{L^2(\R^2)}=1\big\}$$
which we again know exists. 
Note that by the definition of the `normalized' sequence and the estimate given in \eqref{funcspaceupbound} we have that \(\tilde{u}_j\) is a bounded sequence in \(H^1(\R^2)\), thus it admits a subsequence \(\tilde{u}_{j_\ell} \) that weakly converges to some $\tilde{u} \in H^1(\R^2)$. On the other hand, the restriction to $D$ of $\tilde{u}_{j_\ell}$ is in $H^1(D)$, so by the Rellich compactness theorem (up to extracting a further subsequence)
$$\tilde{u}_{j_\ell} \to \tilde{u} \,\,  \text{ as $\ell \to \infty$  in  $L^2(D)$}  \,\, \text{ which gives that } \quad \| \tilde{u}_{j_\ell}\|_{L^2(D)}^2 \to \| \tilde{u}\|_{L^2(D)}^2 \,\,  \text{ as $\ell \to \infty$}. $$
Moreover, the weak convergence implies
$$\| \tilde{u}\|_{L^2(\R^2)}^2 \leq \liminf\limits_{\ell\to\infty} \| \tilde{u}_{j_\ell}\|_{L^2(\R^2)}^2 = \liminf\limits_{\ell\to\infty} \left(2\| \tilde{u}_{j_\ell}\|_{L^2(D)}^2 - 1\right) = 2 \| \tilde{u}\|_{L^2(D)}^2 - 1.$$ 
So we have that $2\| \tilde{u}\|^2_{L^2(D)}- \| \tilde{u} \|^2_{L^2(\R^2)} \geq 1$. Therefore, the limit $\tilde{u}$ satisfies the constraint and we have that  
$$\|\nabla \tilde{u}\|_{L^2(\R^2)}^2 \geq \mathcal{R}_D [\tilde{u}] \geq   \inf\limits_{\substack{u\in H^1(\R^2)}} \left\{  \mathcal{R}_D  [u] \,\, : \,\, 2\| u \|^2_{L^2(D)}- \| u \|^2_{L^2(\R^2)}>0\right \}.$$
By the construction of the sequence \(\tilde{u}_{j_\ell}\), we see that 
\begin{align*} 
\mathcal{R}_D [\tilde{u}] \leq \|\nabla \tilde{u}\|_{L^2(\R^2)}^2 & \leq \liminf\limits_{\ell\to\infty} \|\nabla \tilde{u}_{j_\ell}\|_{L^2(\R^2)}^2 \\
& =  \liminf\limits_{\ell \to \infty} \mathcal{R}_D [u_{j_\ell}] \\
& = \inf\limits_{\substack{u \in H^1(\R^2)}} \left\{  \mathcal{R}_D  [u] \,\, : \,\, 2\| u \|^2_{L^2(D)}- \| u \|^2_{L^2(\R^2)}>0  \right\}.
\end{align*}
Therefore, \(\tilde{u}\) is a minimizer of \eqref{eqn_inf_quotient}. 

It only remains to show that the minimizer of \eqref{eqn_inf_quotient} satisfies equation \eqref{eqn_trans_eig_transformed}. Without loss of generality, we assume that the minimizer is `normalized' such that 
$$ 2\| {u} \|^2_{L^2(D)}- \| {u} \|^2_{L^2(\R^2)}=1$$
Denote \(\tau\) as the Lagrange multiplier, which implies that the Lagrangian of our constrained minimization problem is given by 
$$(u,\tau) \longmapsto \|\nabla {u}\|_{L^2(\R^2)}^2 - \tau \left( 2\| u \|^2_{L^2(D)}- \| u \|^2_{L^2(\R^2)}-1 \right)$$
Taking the first variation with respect to $u$ and setting it equal to zero at the minimizer we obtain 
$$\int_{\R^2} \nabla {u}\cdot \nabla\phi \dd x - \tau \int_{\R^2} (2\cdot\mathbf{1}_D-1){u}\phi \dd x = 0 \quad \text{for all} \,\,  \phi \in H^1(\R^2).$$
This is the variational formulation given in \eqref{varform-4RQ} that is equivalent to the clamped transmission eigenvalue problem. So the minimizer $u$ is an eigenfunction and $\tau=k^2$ is the corresponding eigenvalue. Letting $\phi =  {u}$ (eigenfunction) in the variational form gives that 
$$k^2 = \mathcal{R}_D [u]  \quad \text{ which is minimized when} \quad  k_1^2 = \mathcal{R}_D [u_1]$$
where $(k_1 , u_1)$ is the first clamped transmission eigenpair. 
\end{proof}

With our new formulation of the clamped transmission eigenvalue problem given by equation \eqref{eqn_trans_eig_transformed}, we have shown this classical result that the first eigenpair minimizes the associated Rayleigh quotient. This fact is true for many eigenvalue problems for differential operators, see \cite{p-laplace} for the analysis of an eigenvalue problem for the $p$--Laplacian. To extend this result we work in the function space $H^1(\R^2)$ which is not compactly embedded into $L^2(\R^2)$ while the compact embedding is very useful when studying eigenvalue problems associated with a differential operator on bounded regions. As we see in the previous result, the Gagliardo–Nirenberg (Ladyzhenskaya) interpolation inequality can be used to circumvent the lack of a compact embedding.


\section{Analysis of the Principal Transmission Eigenpair}\label{eigpair-analysis}
Now, in this section we will consider proving some of the analytical properties of the principal clamped transmission eigenpair that was noted in earlier studies \cite{HLK-TEP,HarrisKleefeldLee2025}. In these papers, the authors used different variational formulations to study this eigenvalue problem. Also, in their numerical experiments, the authors noticed a few interesting aspects of the eigenpairs. Here, we wish to show that these numerical observations are true using the Rayleigh quotient $\mathcal{R}_D  [u]$ defined in \eqref{eqn_frac_form} as well as Theorem \ref{thm_first_eig_in_var_form}.

To begin, we first note that Theorem \ref{thm_first_eig_in_var_form} gives a standard scaling result for the first clamped transmission eigenvalue. The scaling result is the initial step in studying the dependence of the clamped transmission eigenvalues with respect to the region $D$. Indeed, for a fixed $D \subset \R^2$ and any $c>0$ we define 
$$c D = \{ y \in \R^2 \,\, : \,\, y=c x \text{ for some $x \in D$} \}.$$
Therefore, by a simple change of variables we can prove the following scaling result. 

\begin{theorem}\label{eig-scaling}
Let $k_1(D)$ and $k_1(c D)$ denote the first clamped transmission eigenvalue of \eqref{eqn_trans_eig_transformed} or equivalently \eqref{tep1}--\eqref{tep2}, corresponding to the region $D$ and the scaled region $c D$, respectively. Then we have that 
\[
k_1(c D) =\frac{1}{c} k_1(D) \quad \text{for any $c>0$.}
\]
\end{theorem}
\begin{proof}
To prove the claim, notice that by Theorem \ref{thm_first_eig_in_var_form} and the change of variables $y=cx$ implies that $\dd y=c^2 \dd x$ and $\grad_y = c^{-1} \grad_x$. With this, we see that the Rayleigh quotient gives
\begin{align*}
k_1^2(c D) \,\, &= \min\limits_{\substack{u\in H^1(\R^2)} } \left\{  \mathcal{R}_{c D}  [u]  \,\, : \,\, 2\| u \|^2_{L^2(c D)}- \| u \|^2_{L^2(\R^2)}>0 \right\}\\
			   &= \min\limits_{\substack{u\in H^1(\R^2)} } \left\{  \frac{1}{c^2} \mathcal{R}_{ D}  [u]  \,\, : \,\, 2\| u \|^2_{L^2( D)}- \| u \|^2_{L^2(\R^2)}>0 \right\}\\
			   &=\frac{1}{c^2} k_1^2(D).
\end{align*}
These calculations give the  result. 
\end{proof}

With this, we can further investigate the dependence of the principal eigenvalue with respect to the region $D$. It was conjectured that the principal eigenvalue $k_1(D)$ is monotone decreasing with respect to the region $D$. This was observed via numerical experiments in both \cite{HLK-TEP,HarrisKleefeldLee2025}. Note that, Theorem \ref{eig-scaling} confirms this when the region is scaled by $c>1$. Since we know that the far--field data can recover the clamped transmission eigenvalue by Theorem 4.1 in \cite{HLK-TEP} this would give an inverse spectral result, that the first eigenvalue determines the area as is case for the first Dirichlet eigenvalue of $- \Delta$ in $D$. With the results from the previous section we can prove the following monotonicity result. 

\begin{theorem}\label{eig-monotone}
Let $k_1(D)$ denote the first clamped transmission eigenvalue of \eqref{eqn_trans_eig_transformed} or equivalently \eqref{tep1}--\eqref{tep2}, corresponding to the region $D$. Then we have that for two regions $D_1 \subseteq D_2$, their first eigenvalues satisfy 
\[
k_1(D_1) \geq k_1(D_2) \quad \text{i.e. monotone decreasing with respect to $D$.}
\]
Moreover, if the set inclusion is strict then so is the inequality for the eigenvalues, provided the regions satisfy that $|D_2 \setminus \overline{D_1}|>0$. 
\end{theorem}
\begin{proof}
This result is a consequence of Theorem \ref{thm_first_eig_in_var_form}. Indeed, we let $u_j$ be the corresponding principal eigenfunctions for the region $D_j$ for $j=1,2$ and  we have that 
\begin{align*}
k^2_1(D_1) = \mathcal{R}_{D_1} [u_1] \geq \mathcal{R}_{D_2} [u_1] \quad \text{since $D_1 \subseteq D_2$.}
\end{align*}
Notice that we have used the fact that $D_1 \subseteq D_2$ and the constraint to obtain 
$$\| u_1 \|^2_{L^2(\R^2)} <  2\| u_1 \|^2_{L^2(D_1)}  \leq 2 \| u_1 \|^2_{L^2(D_2)} \quad \text{which implies} \quad 2\| u_1 \|^2_{L^2(D_2)}- \| u_1 \|^2_{L^2(\R^2)}>0.$$
Continuing to estimate gives that 
\begin{align*}
k^2_1(D_1)&\geq  \min\limits_{\substack{u \in H^1(\R^2) }}\left\{  \mathcal{R}_{D_2} [u] \,\, : \,\,  2\| u \|^2_{L^2(D_2)}- \| u \|^2_{L^2(\R^2)}>0 \right\} =\mathcal{R}_{D_2} [u_2]= k^2_1(D_2).
\end{align*}
This proves the claim. 
\end{proof}

In Theorem \ref{eig-monotone} we have verified the conjecture that the first clamped transmission eigenvalue is monotone decreasing with respect to $D$. This fact was observed numerically in Table 2 in \cite{HLK-TEP} as well as Tables 2 and 4 in \cite{HarrisKleefeldLee2025}. It should also be noted that in the aforementioned manuscripts, the authors noticed this monotonicity for all of the computed eigenvalues. From the monotonicity result given in Theorem \ref{eig-monotone}, we can ask other questions about how the eigenvalues depend on the region $D$. In particular, what happens in the limit as $|D| \to 0^+$. For the case when $D=B_R$ i.e. the ball of radius $R>0$ centered at the origin, we can use separation of variables on \eqref{tep1}--\eqref{tep2} to show that $k>0$ is a clamped transmission eigenvalue if and only if for some $\ell \in \N \cup \{ 0 \}$
\begin{align*}
\det{\begin{pmatrix}
J_{\ell}(kR) & H_{\ell}^{(1)}(\mathrm{i}kR) \\
kJ'_{\ell}(kR) & \mathrm{i}k {H_{\ell}^{(1)}}'(\mathrm{i}kR)
\end{pmatrix}
} = \mathrm{i}k J_{\ell}(kR){H_{\ell}^{(1)}}'(\mathrm{i}kR) - k H_{\ell}^{(1)}(\mathrm{i}kR)J'_{\ell}(kR) = 0
\end{align*}
where $J_{\ell}$ is the Bessel function of first kind with order ${\ell}$ and $H_{\ell}^{(1)}$ is the Hankel function of first kind with order ${\ell}$. With this, we see that the clamped transmission eigenvalues satisfy
\begin{align}\label{firsteig-disk}
k(B_R) = \frac{k(B_1)}{R} \quad \text{which implies} \quad k^2(B_R) = \frac{\pi k^2(B_1)}{\pi R^2} = \frac{\pi k^2(B_1)}{|B_R|}
\end{align}
for any $R>0$. Therefore, for this simple case we see that by \eqref{firsteig-disk} as $|B_R| \to 0^+$ the clamped transmission eigenvalues $k(B_R)\to \infty$. We now prove this limiting result for the clamped transmission eigenvalues for a generic region $D$.

\begin{theorem}\label{eig-smallvolumelimit}
Let $k(D)$ denote a clamped transmission eigenvalue of \eqref{eqn_trans_eig_transformed} or equivalently \eqref{tep1}--\eqref{tep2}, corresponding to the region $D$. Then we have that as 
\[ |D| \to 0^+ \quad \text{the eigenvalues satisfy } \quad k(D)\to \infty.\]
Moreover, for $k_1(D)$ denoting the first clamped transmission eigenvalue
\[k^2_1(D) \geq \frac{C}{|D|} \quad \text{for some constant $C>0$.}\]
\end{theorem}
\begin{proof}
To prove the claim, we only need to prove the given estimate for the first clamped transmission eigenvalue. Recall that from the proof of Theorem \ref{thm_first_eig_in_var_form} we have the lower bound of the Rayleigh quotient given by 
$$\mathcal{R}_D [u] \geq \frac{C}{ |D|} \quad \text{for all $u \in H^1(\R^2)$ : } \,\, 2\| u \|^2_{L^2(D)}- \| u \|^2_{L^2(\R^2)}>0$$
which gives that 
$$k_1^2 (D) = \mathcal{R}_D [u_1]\quad \text{ and therefore } \quad k_1^2 (D) \geq  \frac{C}{ |D|}$$
where $u_1$ is the  corresponding eigenfunction for $k_1$, proving the claim.  
\end{proof}

Notice that the inequality given in Theorem \ref{eig-smallvolumelimit} gives an analog to the explicit result that we have for the case when $D$ is a ball. Theorems \ref{eig-scaling}--\ref{eig-smallvolumelimit} give analytical results on how the clamped transmission eigenvalues depend on the region $D$. In particular, the lower bound with respect to the measure of $D$ given in Theorem \ref{eig-smallvolumelimit} is sharp for the ball as our discussion has shown. 

Another question is if the estimate is sharp for a generic region for some appropriate assumptions on $D$. To this end, we first define the interior--radius of the region $D$
$$ \rho_D = \sup\left\{ \rho>0 \,\, : \,\, \exists x \in D \, \text{ where } \, B(x ; \rho) \subset D \right\}.$$
Here, we let $B(x ; \rho)$ denote the ball centered at $x$ with radius $\rho>0$.  For any non--empty, simply connected bounded open set $D$, we have that 
$$\pi \rho_D^2 \geq \alpha |D| \quad \text{ for some  $0<\alpha \leq1$}.$$
Note, that in general the value of $\alpha$ will depend on the set $D$. In particular, for the `flat' ellipse whose boundary is given by 
$$\partial D = (\cos t , b \sin t)^\top  \quad \text{for $0<b <1$ and $t \in [0,2\pi)$}$$ 
then it is clear that for the `flat' ellipse $\alpha = b$. With this in mind, for this region we can prove an upper bound similar to the lower bound given in the previous result. To this end, we can choose $r>0$ such that  
$$\rho_D/2 < r <  \rho_D \quad \text{ which give that } \quad \exists \, B(x_0 ; r) \subset D.$$ 
Now, take $\psi \in C^\infty_c \left(B_1 \right)$ and by extending by zero outside of $B_1$ define 
$$\varphi (x) = \psi \left( \frac{x-x_0}{r} \right) \in H^1(\R^2) \quad \text{satisfying} \quad  2 \| \varphi \|^2_{L^2(D)}- \| \varphi \|^2_{L^2(\R^2)}>0.$$ 
with supp$(\varphi) \subseteq B(x_0 ; r) \subset D$. From direct calculations we have that 
$$\| \varphi \|^2_{L^2(\R^2)}=\| \varphi \|^2_{L^2(D)} = r^2 \| \psi \|^2_{L^2(B_1)}  \quad \text{and} \quad \| \grad \varphi \|^2_{L^2(\R^2)} =  \| \grad \psi \|^2_{L^2(B_1)}$$ 
via a simple change of variables. Using the estimate for $r$ and the interior--radius, we can see that 
$$ r^2> \rho^2_D/4 \geq \frac{\alpha}{4\pi} |D|.$$
Theorem \ref{thm_first_eig_in_var_form} gives that 
\begin{align} \label{upboundD1}
k^2_1(D) \leq \mathcal{R}_D[\varphi] &= \frac{\| \grad \varphi \|^2_{L^2(\R^2)}}{\| \varphi \|^2_{L^2(D)} } =\frac{\| \grad \psi \|^2_{L^2(B_1)}}{r^2 \| \psi \|^2_{L^2(B_1)} }\leq \frac{4\pi}{\alpha |D|} \frac{\| \grad \psi \|^2_{L^2(B_1)}}{ \| \psi \|^2_{L^2(B_1)} } =  \frac{C}{\alpha |D|}
\end{align} 
where the constant $C>0$ is independent of $D$. We recall that $\alpha$ in the above estimate does depend on $D$. For the `flat' ellipse, since $\alpha=b$ we see that for this region get the estimate 
$$k^2_1(D) \leq \frac{C}{|D|^2} \,\,\, \text{ since $|D|=\pi b$ for the `flat' ellipse.} $$

This upper bound is a product of the fact that the region is `flat'. For this example, the diameter of $D$ is fixed but the interior--radius tends to zero for $b\ll 1$. This seems to be the degenerate case. If we make an assumption that $D$ satisfies the shape--regularity condition (standard in studying FEMs \cite{SunZhouFEM}) that for the diam$(D)=d_D$ we let
$$D \in \Sigma_\beta :=\left\{ D \subset \R^2 \,\, : \,\, \text{for a fixed constant $\beta>0$} \,\, \text{we have that} \,\,  \frac{d_D}{\rho_D} \leq \beta \right\}$$
then for any $D \in \Sigma_\beta$ the upper bound in \eqref{upboundD1} becomes 
\begin{align} \label{upboundD2}
k^2_1(D) \leq \frac{\| \grad \psi \|^2_{L^2(B_1)}}{r^2 \| \psi \|^2_{L^2(B_1)} } \, \leq \, \frac{\pi \beta^2}{ |D|} \frac{\| \grad \psi \|^2_{L^2(B_1)}}{ \| \psi \|^2_{L^2(B_1)} } =  \frac{C_\beta}{|D|}
\end{align} 
where we have used the estimate $|D| \leq (\pi/4)d_D^2$ and $\beta$ is a universal constant for all $D$ in the set $\Sigma_\beta$. In the following section, we will numerically investigate the rate of convergence i.e. study if $k^2_1(D) = \mathcal{O}(1/|D|^p)$ as $|D| \to 0^+$ for some $p \geq 1$. We will consider the convergence for a `flat' region as well as shape regular regions $D \in \Sigma_\beta $ such that $|D| \to 0^+$. 

Now, we note that Theorem 3.8 in \cite{HarrisKleefeldLee2025} proves that for $k_1(D)$ denoting the first clamped transmission eigenvalue of $D$ and $\lambda_1(D)$ denoting the first Dirichlet eigenvalue of $- \Delta$ in $D$ we have that  $k^2_1(D) \leq \lambda_1(D)$. For other transmission eigenvalue problems coming from other scattering phenomena  (see, e.g., \cite{Fiora-TEpaper,Rafa-TEpaper,Victor-TEpaper,Xiang-TEpaper}) lower bounds for the first transmission eigenvalue are proven with respect to the first Dirichlet eigenvalue. Connecting different eigenvalues for a region $D$ is a classical area of study. Recently, in \cite{MagneticLaplacian-EigInequ} an interlacing inequality for the Dirichlet and Neumann eigenvalues of the magnetic Laplacian was proven. We will now use the Rayleigh quotient to show that the first clamped transmission eigenvalue is strictly less than the first Dirichlet eigenvalue.

\begin{theorem}\label{TE-lessthanDirEig}
Let $k_1(D)$ denote the first clamped transmission eigenvalue of \eqref{eqn_trans_eig_transformed} or equivalently \eqref{tep1}--\eqref{tep2}, corresponding to the region $D$. Then we have that
\[k^2_1(D) < \lambda_1(D)\]
where $\lambda_1(D)$ is the first Dirichlet eigenvalue of $- \Delta$ in $D$.
\end{theorem}
\begin{proof}
In order to prove the claim, we will show that $\lambda_1(D)$ is not a minimal value of the Rayleigh quotient. To this end, denote \(\phi_{\text{Dir}}\in H^1_0(D)\) as the corresponding Dirichlet eigenfunction. It is well known that the first Dirichlet eigenfunction has a single sign in $D$ and is non--trivial. Therefore, by scaling we can assume that 
$$\phi_{\text{Dir}}> 0 \quad \text{a.e. in $D$ such that } \,\, \int_D \phi_{\text{Dir}}(x)\dd x = 1.  $$
Due to the zero trace on $\partial D$, we can extend $\phi_{\text{Dir}}$ by zero into $\exD$. 

Now, take any smooth, non--negative cutoff function \(\chi\in C^\infty_c(\R^2)\) such that \(\chi(x) = 1\) for any \(x\in D\).
Now, we consider the perturbation of the extended Dirichlet eigenfunction
$${\varphi}_t (x):= \phi (x) + t\chi(x) \in H^1(\R^2),\qquad \text{ for any } \,\,\, t \in \R$$
where $\phi (x) = \phi_{\text{Dir}}(x)$ for $x \in D$ and $\phi(x) = 0$ for $x \in \exD$. By direct calculations, we see that the constraint  
$$  2\|{\varphi}_{t} \|_{L^2(D)}^2 - \|{\varphi}_{t} \|_{L^2(\R^2)}^2 = \|\phi_{\text{Dir}}\|_{L^2(D)}^2 + 2t +\left( |D|-\| \chi \|^2_{L^2(\exD)}\right) t^2.$$
This is due to the fact that 
$${\varphi}_{t}(x) = \begin{cases}
            \phi_{\text{Dir}}(x)+t, & x\in D, \\
            t \chi(x), & x\in \exD .
        \end{cases} $$
Notice that at $t=0$ we have $ 2\|{\varphi}_{0} \|_{L^2(D)}^2 - \|{\varphi}_{0} \|_{L^2(\R^2)}^2 = \|\phi_{\text{Dir}}\|_{L^2(D)}^2 >0.$         
Therefore, since the above equality for the constraint of ${\varphi}_{t}$ is continuous for all $t$ we have that there exists $\tau>0$ such that $2\|{\varphi}_{t} \|_{L^2(D)}^2 - \|{\varphi}_{t} \|_{L^2(\R^2)}^2 >0$ for all $t \in [-\tau,\tau]$.

Using the fact that 
$$ \nabla{\varphi}_{t}(x) = \begin{cases}
            \nabla \phi_{\text{Dir}}(x), & x\in D, \\
            t\nabla\chi(x), & x\in \exD
        \end{cases}$$
we can directly compute
$$ \mathcal{R}_D [{\varphi}_{t}]  = \frac{\|\nabla \phi_{\text{Dir}}\|_{L^2(D)}^2 + t^2\|\nabla\chi\|_{L^2(\exD)}^2}{\|\phi_{\text{Dir}}\|_{L^2(D)}^2 + 2t +\left( |D|-\| \chi \|^2_{L^2(\exD)}\right) t^2}.$$
 Differentiating with respect to \(t\) and evaluating at \(t=0\), we obtain
$$ \left(\frac{\dd}{\dd t} \mathcal{R}_D [{\varphi}_{t}] \right)\Bigg|_{t=0} = -\frac{2\|\nabla \phi_{\text{Dir}}\|_{L^2(D)}^2}{\|\phi_{\text{Dir}}\|_{L^2(D)}^4} < 0.$$
This implies that ${\varphi}_{0}$ is not a minimizer for the Rayleigh quotient and therefore gives that  
$$k_1^2(D)< \mathcal{R}_D [{\varphi}_{0}] =  \frac{\|\nabla \phi_{\text{Dir}}\|_{L^2(D)}^2}{\|\phi_{\text{Dir}}\|_{L^2(D)}^2} = \lambda_1(D)$$
which gives the result.
\end{proof}

This gives a stronger result when comparing the first clamped and first Dirichlet eigenvalues of $D$. Note that the strict inequality in Theorem \ref{TE-lessthanDirEig} and the monotonicity result given in Theorem \ref{eig-monotone} come from studying the Rayleigh quotient defined in equation \eqref{eqn_frac_form}. Note, that the Rayleigh quotient associated with the variational forms studied in \cite{HLK-TEP,HarrisKleefeldLee2025} are not amenable to proving these results. Also note that, the strict inequality was observed in Tables 1, 2 and 4 of \cite{HarrisKleefeldLee2025} but a non--strict inequality was proven.

Now that we have studied the properties of the principal eigenvalues, we can now study the corresponding eigenfunctions. Study of the analytical properties of the eigenfunctions for other transmission eigenvalue problems can be found in \cite{eigfuncstudy1,eigfuncstudy2,eigfuncstudy3}. In \cite{HLK-TEP}, plots of the real and/or imaginary parts for the computed eigenfunctions $(v,w)$ where presented. One thing to note is that for the three regions that are considered in \cite{HLK-TEP}, we see that the real/imaginary parts of the principal eigenfunctions have a single sign in there respective regions i.e. 
$$ v >0 \text{ (or $<0$) in $D$} \quad \text{ and } \quad  w <0 \text{ (or $>0$) in $\exD$.}$$
From the definition of the combined (real--valued) eigenfunction $u$ defined in \eqref{new-eigfunc} this would translate to $u>0$ (or $<0$) in $\R^2$. The single sign property is well known for the first Dirichlet and Neumann eigenfunction of $- \Delta$ in $D$. We now prove a similar result for the principal clamped transmission eigenfunction.

\begin{theorem}\label{eigfunc-sign}
Let $u_1$ denote the eigenfunction corresponding to the first clamped transmission eigenvalue of \eqref{eqn_trans_eig_transformed} or equivalently \eqref{tep1}--\eqref{tep2}, corresponding to the region $D$. If we assume that $u_1> 0$ (or $< 0$) on $\partial D$ then we have that $u_1 > 0$ (or $< 0$) in $\R^2$.
\end{theorem}
\begin{proof}
To prove the claim, we will assume that $u_1> 0$ on $\partial D$ and the other case can be proven similarly. Now, with this assumption, we have that 
$$ -\Delta u_1 = k_1^2 u_1 \,\,\text{ in } D \quad \text{and } \quad u_1  > 0\,\, \text{on } \partial D.$$
With this, we have that by elliptic regularity (see, e.g., Chapter 8 in \cite{salsa}), $u_1$ is continuous in the interior of $D$. We also have that $u^* = \min\{u_1 , 0\} \in H^1_0(D)$ due to the sign assumption on the boundary $\partial D$. Therefore, we have that 
$$\int_D \grad u_1 \cdot \grad u^* \, \dd x = k^2_1\int_D  u_1   u^* \, \dd x \implies \int_{\text{supp}(u^* )} \grad u_1 \cdot \grad u^* \, \dd x = k^2_1\int_{\text{supp}(u^* )}  u_1   u^* \, \dd x. $$
We can use the above equality and Theorem \ref{TE-lessthanDirEig} to prove the result. With this in mind, notice that $u_1 = u^*$ in  $\text{supp}(u^* ) \subset D$. 
This implies that by the above equalities 
$$ \int_{\text{supp}(u^* )} | \grad u^*|^2 \, \dd x = k^2_1 \int_{ \text{supp}(u^* )}  |u^*|^2  \, \dd x \implies \int_{D} | \grad u^*|^2 \, \dd x = k^2_1 \int_{D}  |u^*|^2  \, \dd x.$$ 
By the Poincaré inequality for $H^1_0(D)$
$$ \lambda_1 \int_{ D}  |u^*|^2  \, \dd x \leq \int_{D} | \grad u^*|^2 \, \dd x =  k^2_1 \int_{D}  |u^*|^2  \, \dd x$$ 
which gives that $(\lambda_1 - k^2_1) \|u^*\|^2_{L^2(D)} \leq 0$. By Theorem \ref{TE-lessthanDirEig} we have that  $\|u^*\|^2_{L^2(D)} = 0$ implying that $u^* = 0$ a.e. in $D$ and by the continuity $u_1 >0$ in $D$.  

Now, we focus on $\exD$ and we proceed by way of contradiction. So assume that there exists $\Omega\subset \exD$ such that $u_1 < 0$ in $\Omega$ such that $|\Omega|>0$. Again we have appealed to elliptic regularity to obtain the continuity of $u_1$ in $\exD$ which implies that $\|u_1\|_{L^2(\Omega)}>0$. We now similarly consider
$$u^* = \frac{u_1 + |u_1|}{2}$$
that is, removing the negative part in \(u_1\) and reset it to zero and that almost everywhere
$$\nabla|u_1| = \begin{cases}
            \frac{u_1}{|u_1|}\nabla u_1, & u_1\neq 0, \\
            \,\,\,\, 0, & u_1=0.
        \end{cases}$$
Note that we have that $u^* \in H^1(\R^2)$ and by definition
$$ \|\nabla u^* \|_{L^2(\R^2)} \leq \frac{\|\nabla u_1\|_{L^2(\R^2)} + \|\nabla |u_1|\|_{L^2(\R^2)}}{2} \leq \|\nabla u_1\|_{L^2(\R^2)}.$$
Moreover, since $\|u^* \|_{L^2(\Omega)} = 0 $ we see that $\|u^* \|_{L^2(\Omega)}< \|u_1 \|_{L^2(\Omega)}$,
and by the fact that $|u^*|< |u_1|$, we have $\|u^*\|_{L^2(\R^2)} < \|u_1\|_{L^2(\R^2)}$. Note that we are using the fact that we have already shown $u_1 >0$ in $D$. With these inequalities 
$$ 2\|u^* \|^2_{L^2(D)}-\|u^*\|^2_{L^2(\R^2)}>2\|u_1 \|^2_{L^2(D)}-\|u_1\|^2_{L^2(\R^2)}>0$$
so $u^*$ satisfies the constraint. This gives that
$$\mathcal{R}_D [u^*]= \frac{\|\nabla u^*\|_{L^2(\R^2)}^2}{2\|u^*\|_{L^2(D)}^2 - \|u^*\|_{L^2(\R^2)}^2} < \frac{\|\nabla u_1\|_{L^2(\R^2)}^2}{2\|u_1\|_{L^2(D)}^2 - \|u_1\|_{L^2(\R^2)}^2} = \mathcal{R}_D [u_1].$$
This gives a contradiction, as by Theorem \ref{thm_first_eig_in_var_form} the Rayleigh quotient is minimal at the eigenfunction $u_1$ but the above inequalities implies that $u^*$ gives a smaller value of the Rayleigh quotient. 
\end{proof}

Note that the assumption that the eigenfunction has one sign on the boundary is needed in the above result. Due to the denominator of the Rayleigh quotient the argument in proving Theorem \ref{eigfunc-sign} for the case $u_1>0$ in $\exD$ cannot be applied in all of $\R^2$ if the eigenfunction can change sign across  the boundary. 
With the results in this section, we have provided analytical results for the principal clamped transmission eigenpair. In particular, we investigated the dependence on $D$ where we have shown monotonicity and the limiting behavior as the size tends to zero for the principal eigenvalue. For the principal eigenfunctions we have proven that they do not change sign in their respective regions akin to the results for Dirichlet and Neumann eigenfunctions. 

Now that we have Theorem \ref{eigfunc-sign}, we can assume that the eigenfunction $u_1$ is positive in $\R^2$. This is a key ingredient needed to prove our next result. With the results we have already proven, we can now prove a Faber--Krahn type inequality for the clamped transmission eigenvalues. This inequality gives that the ball minimizes the value of the first clamped transmission eigenvalue over regions with a fixed measure. Similar result for other eigenvalue problems have recently been studied in \cite{BogoselBucur,mainref,p-laplace}. 

\begin{theorem}\label{FK-inequ}
Let $k_1(D)$ and $k_1(B_R)$ denote the first clamped transmission eigenvalue of \eqref{eqn_trans_eig_transformed} or equivalently \eqref{tep1}--\eqref{tep2}, corresponding to the region $D$ and the ball $B_R$, respectively. Assume that $|D|=|B_R|$, then we have the estimate 
\[k_1(B_R) \leq k_1(D)\]
provided that the corresponding eigenfunction for $D$ has a fixed sign on $\partial D$.
\end{theorem}
\begin{proof}
To prove the claim, we argue similarly to the case of the Dirichlet eigenvalue. So for the region $D$ we let $B_R$ be its symmetric rearrangement. Then we can let $u_1 \in H^1(\R^2)$ be the principal eigenfunction associated to $k_1(D)$ (that we can assume is positive) and let $u^\sharp \in H^1(\R^2)$ be its corresponding symmetric--decreasing rearrangement. From this we have the estimates (see, e.g., Chapter 3 in \cite{lieb2001analysis})
$$\|u_1 \|^2_{L^2(D)} \leq \| u^\sharp \|^2_{L^2(B_R)} \quad \text{ and } \quad \|u_1 \|^2_{L^2(\R^2)} = \| u^\sharp \|^2_{L^2(\R^2)}$$ 
which implies that 
$$0< 2\|u_1 \|^2_{L^2(D)} - \|u_1 \|^2_{L^2(\R^2)} \leq  2\| u^\sharp \|^2_{L^2(B_R)} -  \| u^\sharp \|^2_{L^2(\R^2)}.$$ 
This gives that the symmetric--decreasing rearrangement $u^\sharp$ satisfies the constraint for the Rayleigh quotient $\mathcal{R}_{B_R} [ \, \cdot \, ]$. Now, by the appealing to the P\'olya–Szeg\H{o} inequality (see, e.g., Chapter 7 in \cite{lieb2001analysis}) given by  
$$ \|\nabla u^\sharp \|_{L^2(\R^2)}^2 \leq \|\nabla u_1\|^2_{L^2(\R^2)} $$ 
and using the Rayleigh quotient along with Theorem \ref{thm_first_eig_in_var_form} we obtain that 
$$k^2_1(B_R) \leq  \mathcal{R}_{B_R} [ u^\sharp ] \leq  \mathcal{R}_D [u_1] = k^2_1(D)$$ 
which proves the inequality. 
\end{proof}

With the results in these sections, we have provided a plethora of analytical results for the clamped transmission eigenvalues. Here, we have studied the principal eigenpair to show that they have similar properties to the principal Dirichlet eigenpair. The principal eigenvalues are monotone with respect to $D$ and the corresponding eigenfunctions have a single sign. We have also provided a lower bound for the eigenvalues with respect to $|D|$ as well as investigated a possible upper bound for the principal eigenvalue with respect to $|D|$. In the following section, we will provide numerical evidence for these results and consider other possible areas for future study for the clamped transmission eigenvalues.

\section{Numerical Investigation}\label{numerics}
In this section, we will provide numerical examples for the theory developed in the previous sections. In order to compute the clamped transmission eigenvalues and eigenfunctions we will use the method described in Section 5 of \cite{HLK-TEP}. Therefore, we note that this is an eigenvalue problem in an unbounded domain, so we use boundary integral equations to approximate the eigenvalues and eigenfunctions. Here, in order to compute the eigenpairs $k>0$ and $(v,w)  \in H^1(D) \times H^1(\exD) $ we use boundary integral operator where we make a single--layer ansatz. This gives a system of integral equations that depends non--linearly on $k$, where we use the method introduced in \cite{beyn} to compute the eigenpairs.

In our study, we will first consider three shapes for the region $D$. Since we are using boundary integral equations we will take regions that have smooth boundaries in our examples. The regions that we will initially consider are given in Table \ref{regions1}. 
\begin{table}[H]
	\centering
	\begin{tabular}{ l l }
\hline
\hline
		Scatterer     &  Parameterization \\
\hline 
\hline
		\vspace{-2.0ex}\\
		Star & $x(t)={\displaystyle 0.25(3\cos(5t)/10+2) \big(\cos(t),\sin(t) \big)^\top}$
		\vspace{1.5ex} \\ 
		Peanut  & $x(t)={\displaystyle 0.5\sqrt{3\cos^2(t)+1} \big(\cos(t),\sin(t) \big)^\top}$  
		\vspace{1.5ex}\\ 
		Kite & $x(t)={\displaystyle \big(0.75\cos(t)+0.3 \cos(2t),\sin(t) \big)^\top}$
		\vspace{1.5ex}\\
\hline
	\end{tabular}\caption{The boundary parameterizations of $\partial D=x(t)$ for $t\in [0,2\pi)$.} \label{regions1}
\end{table}

We first consider some numerical examples of Theorem \ref{eig-scaling} and \ref{eig-monotone}. These are the result that provided the theoretical justification for the monotonicity of the first eigenvalue with respect to $D$. Theorem \ref{eig-scaling} in particular gives a scaling relationship for the first eigenvalue and Theorem \ref{eig-monotone} give that the first eigenvalue is decreasing with respect to $D$. For the regions given in Table \ref{regions1} we will test the scaling and monotonicity of the clamped transmission eigenvalues (CTEs). Therefore, in Table \ref{table1} we provide the first three eigenvalues for the scaled regions $c_j D$ for the increasing scaling parameters $c_j = (1+j/10)$.
\begin{table}[H]
\centering
 \begin{tabular}{r|r|r|r|r|r|r|r}
                & CTE& $c_0=1$ & $c_1=1.1$ & $c_2=1.2$ & $c_3=1.3$ & $c_4=1.4$ & $c_5=1.5$\\
  \hline
               & $k_1$ & 1.91665& 1.74241& 1.59721& 1.47435& 1.36904&1.27777\\
  Kite       & $k_2$ & 3.38373& 3.07611& 2.81977& 2.60287& 2.41695&2.25582\\
               & $k_3$ & 3.75151& 3.41047& 3.12626& 2.88578& 2.67965&2.50101\\
               \hline 
               & $k_1$ & 3.26715& 2.97014& 2.72263& 2.51320& 2.33368&2.17810\\
   Star     & $k_2$ & 6.18638& 5.62399& 5.15532& 4.75876& 4.41885&4.12426\\
               & $k_3$ & 6.18638& 5.62398& 5.15532& 4.75876& 4.41885&4.12426\\
               \hline 
               & $k_1$ & 2.13093& 1.93721& 1.77577& 1.63917& 1.52209&1.42062\\
Peanut    & $k_2$ & 3.41900& 3.10818& 2.84917& 2.63000& 2.44214&2.27933\\
               & $k_3$ & 4.70289& 4.27536& 3.91908& 3.61761& 3.35921&3.13526\\
               \hline 
 \end{tabular}
 \caption{The first three clamped transmission eigenvalues for the scaled kite--shaped, star--shaped, and peanut--shaped regions given in Table \ref{regions1}.} \label{table1}
\end{table}

The results in Table \ref{table1} confirm the scaling result in Theorem \ref{eig-scaling} not just for the first CTE but for the three reported eigenvalues. The scaling result was proved using the fact that the first eigenvalue minimizes the Rayleigh quotient. We do not have a min--max principle for this eigenvalue problem. This is because the denominator of the Rayleigh quotient does not define an appropriate norm on a subset of $H^1(\R^2)$.  We can now use \eqref{firsteig-disk} and the computed eigenvalues in Table \ref{table1} to verify the Faber--Krahn type inequality given in Theorem \ref{FK-inequ}. With this, since we can compute the area for the regions in Table \ref{regions1} we provide the corresponding principal eigenvalue for the associated ball. In Table \ref{FK-table}, we see that the inequality given in Theorem \ref{FK-inequ} is valid for the three regions considered here. 

\begin{table}[!ht]
\centering
 \begin{tabular}{r|c|c|c}
                         & Kite     & Star         & Peanut   \\
\hline
$k_1(D)$  	        &1.91665  & 3.26715 &2.13093 \\
$k_1(B_R)$  	&1.86443 & 3.21126  &2.04238 \\
\hline
 \end{tabular}
 \caption{\label{FK-table}The first CTE for the kite--shaped, star--shaped, and peanut--shaped regions given in Table \ref{regions1} and the first CTE corresponding to the ball $B_R$ with equal measure. }
\end{table}

\begin{figure}[H]
    \centering
    \includegraphics[width=0.425\linewidth]{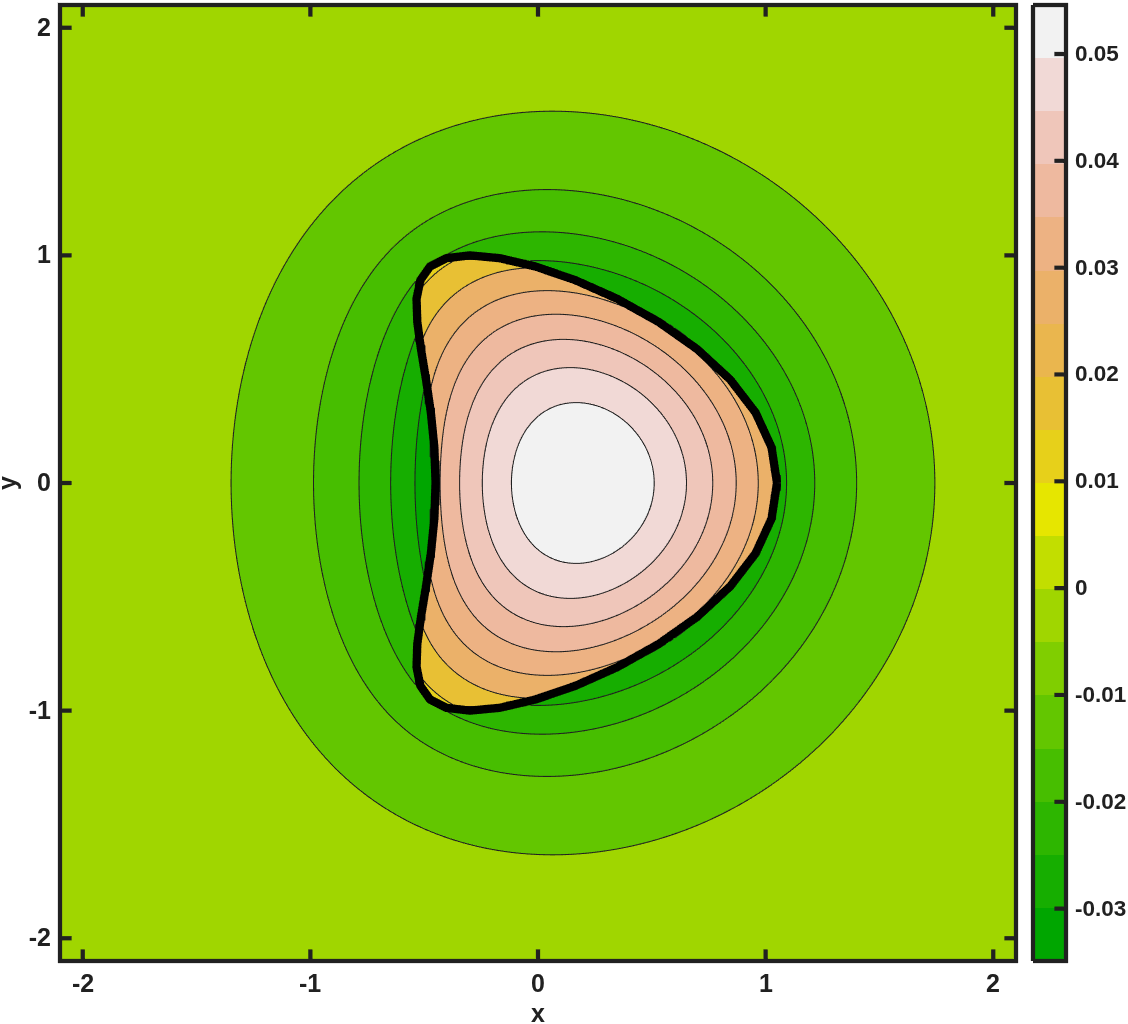}
    \includegraphics[width=0.425\linewidth]{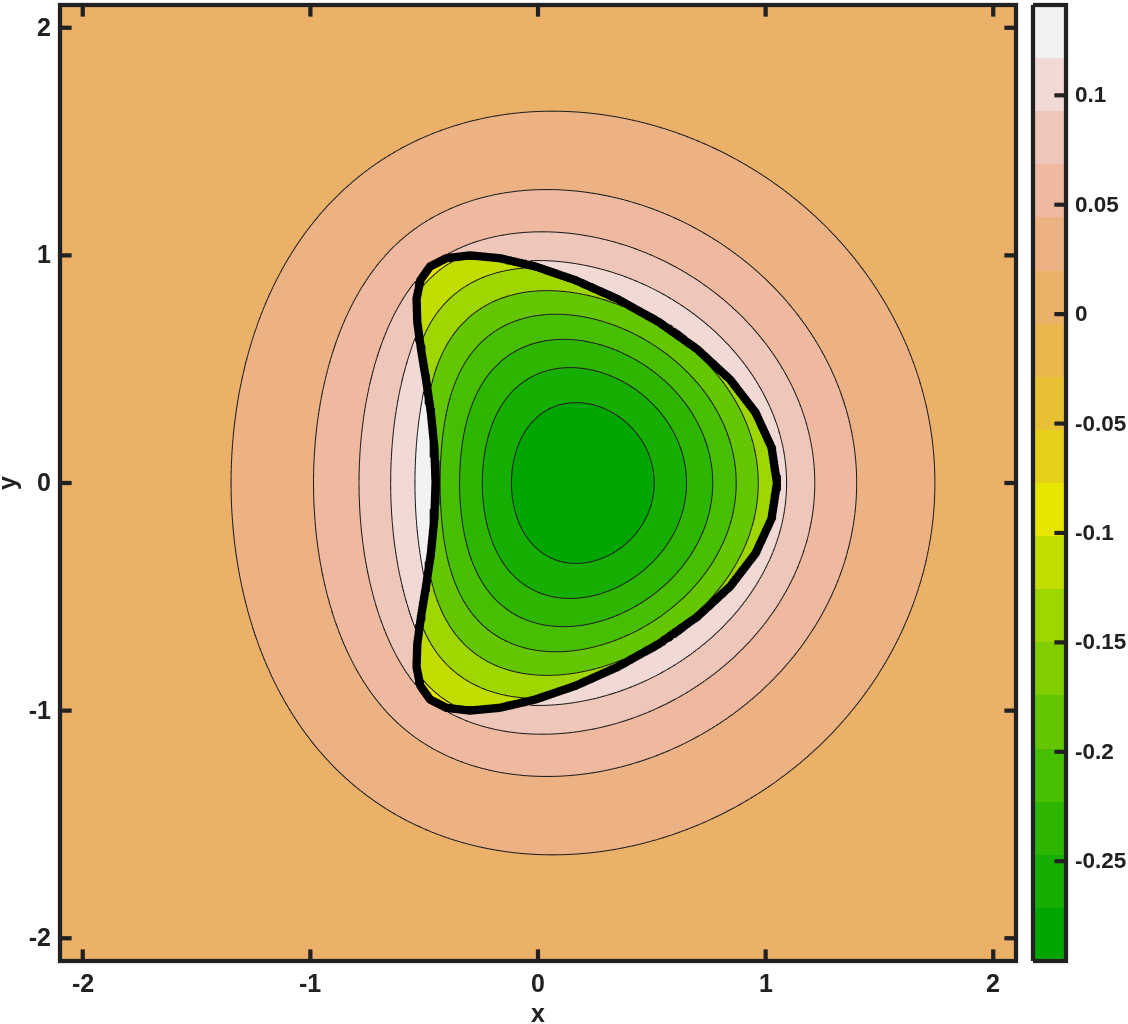}
    \caption{Plots of the real (left) and imaginary (right) parts of the eigenfunctions $w$ within $\mathbb{R}^2\backslash \overline{D}$ and $v$ within $D$ for the kite--shaped region.}
    \label{eigenfctkite}
\end{figure}

Since we use a system of boundary integral equations to compute the CTEs we can also use this to compute the eigenfunctions $(v,w)$. As we have shown in Theorem \ref{eigfunc-sign}, we expect the eigenfunctions corresponding to $k_1$ to have a single sign in their respective domains. More precisely, the real/imaginary parts of $v$ is either positive or negative in $D$ and $w$ is either negative or positive in $\mathbb{R}^2\backslash \overline{D}$. In Figure \ref{eigenfctkite}, we give the plots of the principal eigenfunctions for the kite--shaped region defined in Table \ref{regions1}.

\begin{figure}[H]
    \centering
    \includegraphics[width=0.425\linewidth]{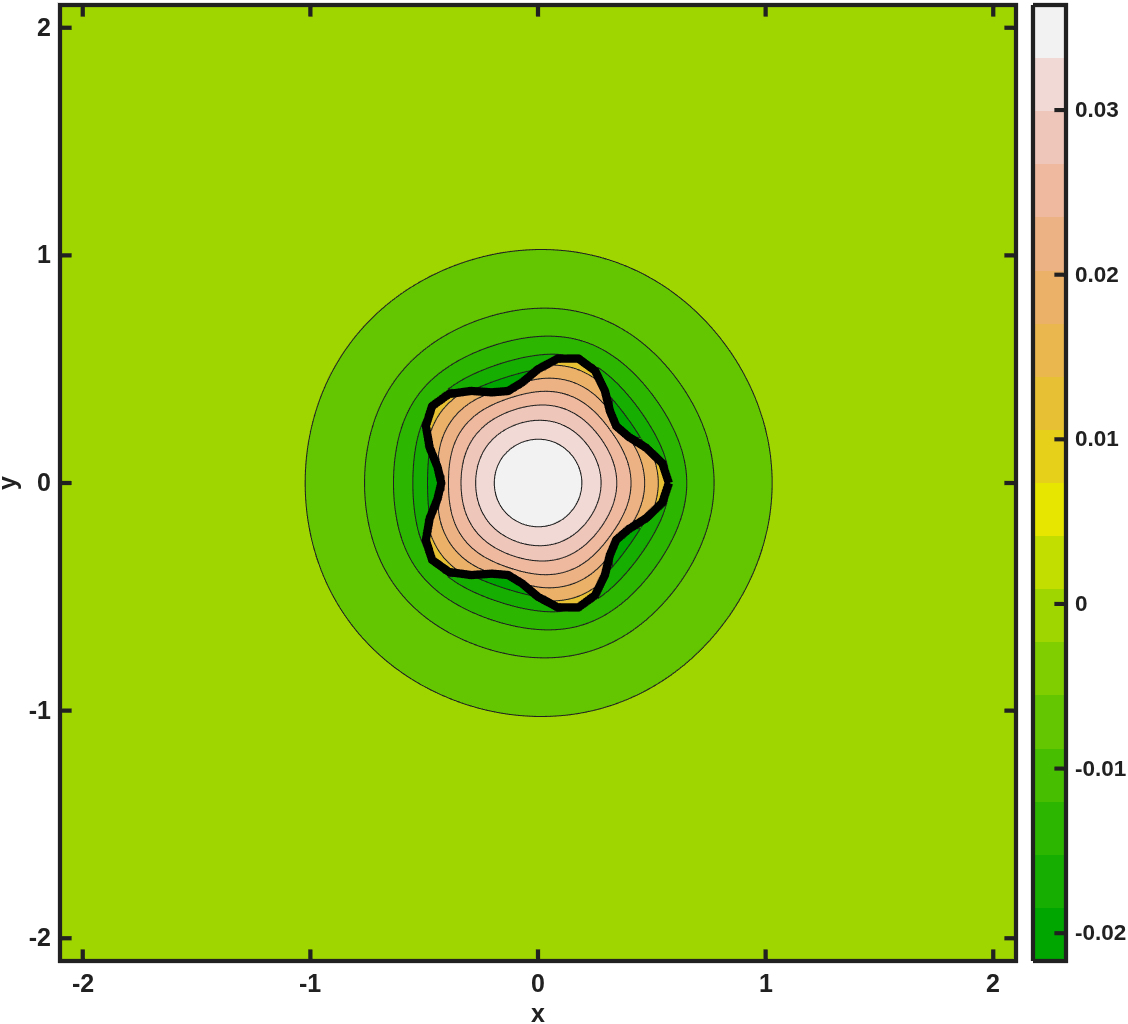}
    \includegraphics[width=0.425\linewidth]{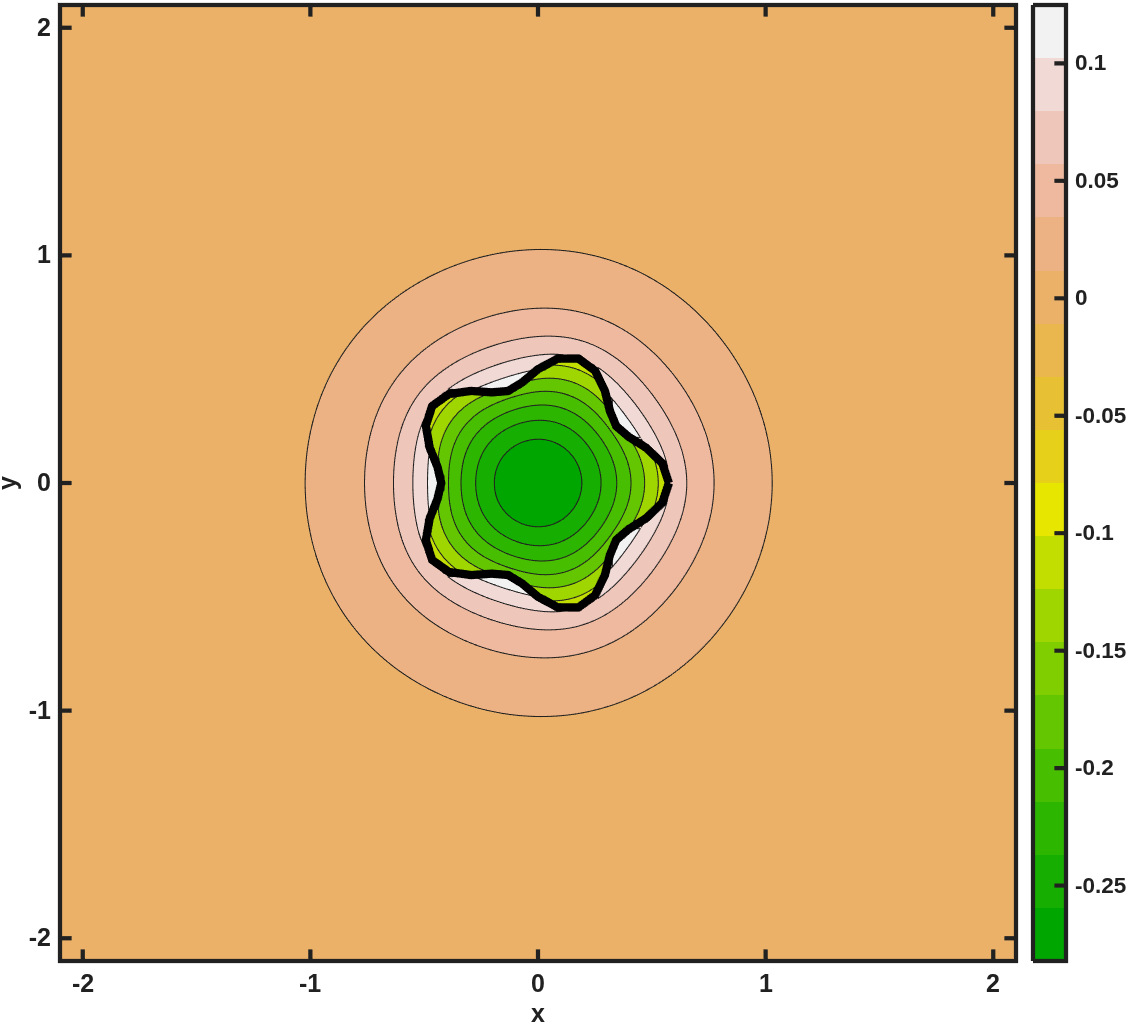}
    \caption{Plots of the real (left) and imaginary (right) parts of the eigenfunctions $w$ within $\mathbb{R}^2\backslash \overline{D}$ and $v$ within $D$ for the star--shaped region.}
    \label{eigenfctstar}
\end{figure}

\begin{figure}[H]
    \centering
    \includegraphics[width=0.425\linewidth]{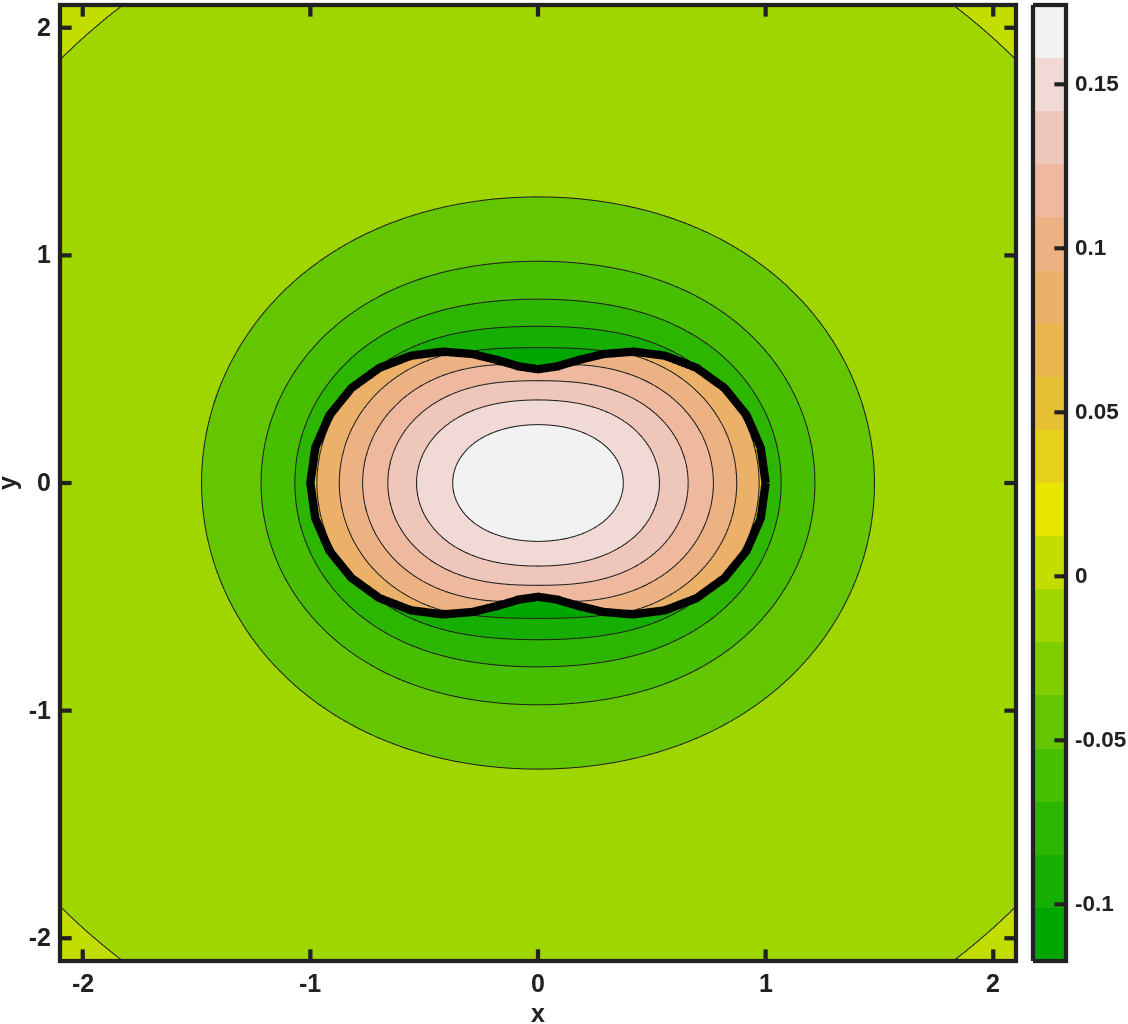}
    \includegraphics[width=0.425\linewidth]{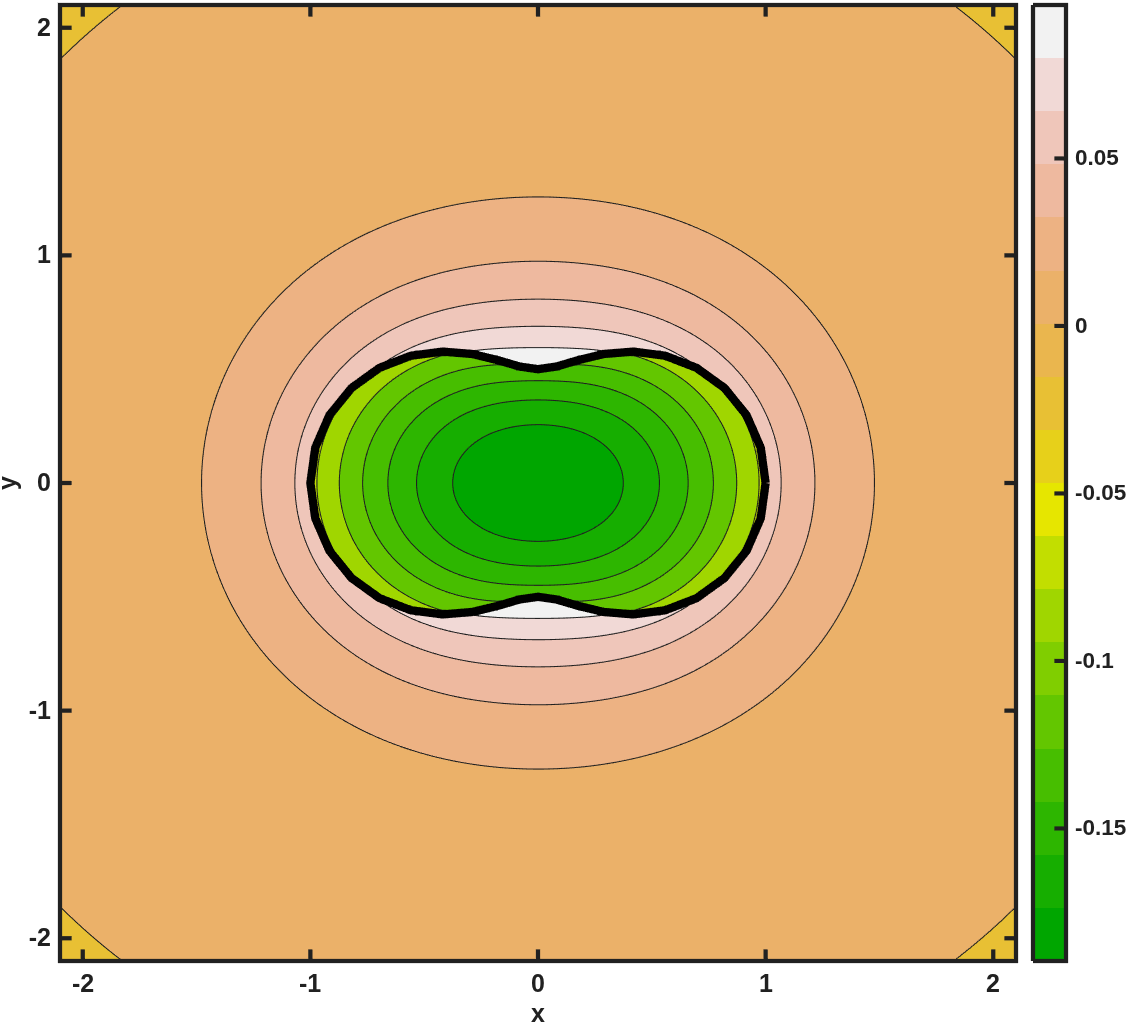}
    \caption{Plots of the real (left) and imaginary (right) parts of the eigenfunctions $w$ within $\mathbb{R}^2\backslash \overline{D}$ and $v$ within $D$ for the peanut--shaped region.}
    \label{eigenfctpeanut}
\end{figure}

We also see the result of Theorem \ref{eigfunc-sign} in full display in Figures \ref{eigenfctstar} and \ref{eigenfctpeanut} for the star--shaped  and peanut--shaped regions whose boundaries are defined in Table \ref{regions1}. Therefore, for each of the three regions we see that the real and imaginary part of the eigenfunctions have a single and opposite sign in their respective regions. In \cite{HLK-TEP}, this is also seen where the real parts of the eigenfunctions for the unit disk and elliptical regions are presented. 


In order to see the dependence of the measure of the set, we give numerical examples for the ellipse that we first considered in Section \ref{eigpair-analysis}. To this end, recall that  
\begin{align} 
\partial D_b = (\cos t , b \sin t)^\top  \quad \text{for $0<b \leq 1$ with $t \in [0,2\pi)$.} \label{regions2}
\end{align}
Now, in Table \ref{table2} we report the first five CTEs for the ellipse. Here, we see that the eigenvalues are monotone increasing with respect to $|D_b|=\pi b$. Note that by Theorem \ref{eig-smallvolumelimit} we know that the CTEs all tend to infinity as $b \to 0^+$.
\begin{table}[H]
\centering
\begin{tabular}{r|c|c|c|c|c|c|c}
CTE & $b=1.0$ & $b=0.9$ & $b=0.8$ & $b=0.7$ & $b=0.6$ & $b=0.5$ & $b=0.4$\\
\hline
$k_1$ &  1.61464 & 1.70401& 1.81492 & 1.95646 & 2.14377 & 2.40418 & 2.79257\\
$k_2$ &  3.05164 & 3.13673& 3.24382 & 3.38233 & 3.56787 & 3.82845 & 4.21987\\
$k_3$ &  3.05164 & 3.30629& 3.62554 & 4.03737 & 4.58853 & 5.26231 & 5.64842\\
$k_4$ &  4.36453 & 4.54060& 4.67571 & 4.82125 & 5.00599 & 5.36324 & 6.53001\\
$k_5$ &  4.36453 & 4.60819& 4.91575 & 5.31515 & 5.85316 & 6.61411 & 7.09087\\
\hline
 \end{tabular}
 \caption{\label{table2}The first five real CTEs for various ellipses with half-axis $(1,b)$ given in \eqref{regions2} using $b\in \{1.0,0.9,0.8,0.7,0.6,0.5,0.4\}$.}
\end{table}
\noindent Notice that we once again see that all of the computed eigenvalues are monotone decreasing with respect to $|D_b|$ that are presented in Table \ref{table2}.

Now we investigate the dependence of the first eigenvalue as $|D| \to 0^+$. We have already shown that all the CTEs will tend to infinity. We are interested in determining the rate at which the (principal) eigenvalue tend to infinity. From Theorem \ref{eig-smallvolumelimit} and the following discussion in Section \ref{eigpair-analysis} we have that there are constants $C_1>0$ and $C_2>0$ such that
\begin{align}\label{RoC-flatellipse}
\frac{C_1}{|D_b|}\leq k_1^2(D_b)\leq \frac{C_2}{|D_b|^2}
\end{align}
for the `flat' ellipse i.e. $0<b<1$. Even though the explicit rate of convergence denoted $p$ was not proven, our analysis would  seem to give that $k^2_1(D_b) = \mathcal{O}(1/|D_b|^p)$ as $|D_b| \to 0^+$ for some $p \in [1,2]$. Now we consider the ellipse such that $b_j = (1-j/10)$ which gives that $|D_j| = |D_{b_j}|$ is decreasing. In order to estimate the rate of convergence  $p$ we define 
\begin{align}\label{RoC-formula}
p_j=-\frac{\log\left(k_1^2(D_{{j}})/k_1^2(D_{{j-1}})\right)}{\log\left(|D_{{j}}|/|D_{{j-1}}|\right)}\,.
\end{align}
This estimation comes from the assumption that $k_1^2(D)\approx (\text{constant})\,|D|^{-p}$ as $|D| \to 0^+$ and solving for the rate of convergence $p$. Hence, the exponent $p$ between two consecutive domains $D_{{j-1}}$ and $D_{{j}}$ can be estimated by the formula in equation \eqref{RoC-formula}. 
\begin{table}[H]
\centering
\begin{tabular}{r|r|r|r|r|r|r}
$b_j$ & $|D_{j}|$ & $k_1(D_{j})$ & $k_1^2(D_{j})$ & $k_1^2(D_{j})\,|D_{j}|$ & $k_1^2(D_{j})\,|D_{j}|^2$ &
rate $p_j$\\
\hline
1.0 & $1.0\pi$ & 1.61464 & 2.607062 & 8.190328 & 25.730674 & ---\\
0.9 & $0.9\pi$ & 1.70401 & 2.903650 & 8.209877 & 23.212881 & 1.0226\\
0.8 & $0.8\pi$ & 1.81492 & 3.293935 & 8.278561 & 20.806292 & 1.0707\\
0.7 & $0.7\pi$ & 1.95646 & 3.827736 & 8.417631 & 18.511336 & 1.1248\\
0.6 & $0.6\pi$ & 2.14377 & 4.595750 & 8.662784 & 16.328964 & 1.1862\\
0.5 & $0.5\pi$ & 2.40418 & 5.780081 & 9.079331 & 14.261779 & 1.2576\\
0.4 & $0.4\pi$ & 2.79257 & 7.798447 & 9.799818 & 12.314814 & 1.3422\\
0.3 & $0.3\pi$ & 3.43775 & 11.818125 & 11.138320 & 10.497619 & 1.4450\\
0.2 & $0.2\pi$ & 4.72953 & 22.368454 & 14.054514 & 8.830712 & 1.5735\\
\hline
\end{tabular}
\caption{Values of first CTE for the `flat' ellipse given by \eqref{regions2} such that $b_j = (1-j/10)$ for $j = 0, \cdots , 8$ and the local rate of convergence $p_j$.}\label{table3}
\end{table}

Notice that the results in Table \ref{table3} verify the analytical rate of convergence in \eqref{RoC-flatellipse} for the  `flat' ellipse. Indeed, for each $j$ we see that the rate of convergence $p_j \in [1,2]$ as expected from the theory. The constants in \eqref{RoC-flatellipse} that are compatible with the data presented in Table \ref{table3} are thus given by
$$C_1=\min_{j = 0, \cdots , 8}  \left\{ k_1^2(D_{j})\,|D_{j}| \right\} =8.1903,$$
and 
$$C_2=\max_{j = 0, \cdots , 8}  \left\{ k_1^2(D_{j})\,|D_{j}|^2\right\} =25.7307$$
giving numerical validation to the theory discussed in the previous sections. 

\begin{figure}[H]
    \centering
    \includegraphics[width=10cm]{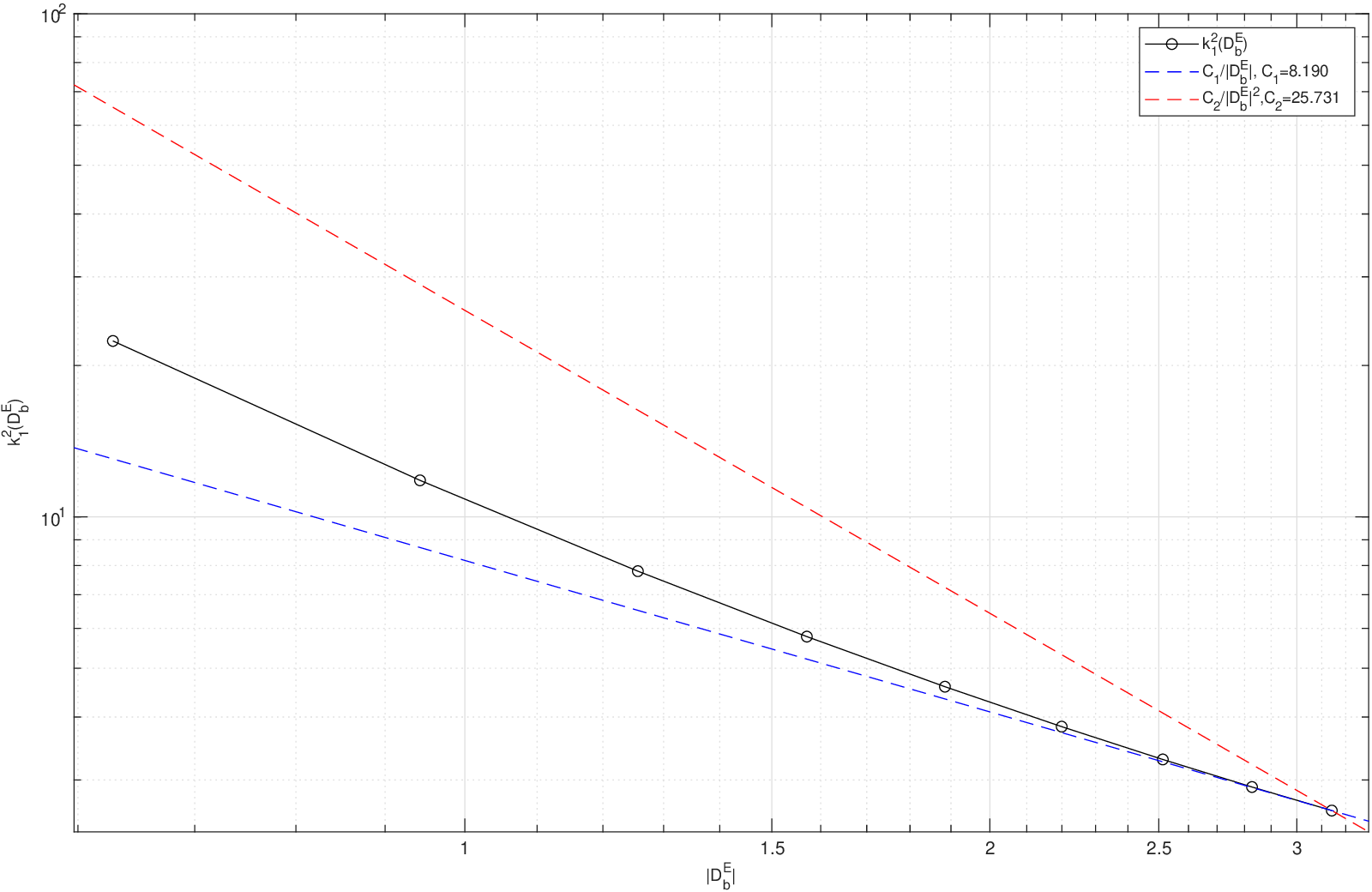}
    \caption{Log--log plot of $k_1^2(D_b)$ versus $|D_b|$ with the theoretical bounds also plotted, i.e.
$C_1/|D_b|$ (line of slope $-1$, $C_1=8.190$) and $C_2/|D_b|^2$ (line of slope $-2$,
$C_2=25.731$).}
    \label{loglog}
\end{figure}

Here we see that the rate of convergence grows from $p_j\approx1.02$ for $j=1$ to $p_j \approx1.57$ for $j=8$, i.e. the data move away from the lower bound towards the upper bound as the area decreases. To visualize this, we provide the log--log plot for the mapping $|D_{b}| \mapsto k_1^2$ in Figure \ref{loglog}. For this extreme case of a thin region, we see that the rate of convergence seems to be bounded above by 2 which may indicate that this is an upper bound for the rate of convergence for a general region giving the upper bound in equation \eqref{RoC-flatellipse}. 

We can do a similar numerical investigation for the scaled regions studied in our previous examples. To this end, we provide the rate of convergence for the scaled kite--shaped, star--shaped and peanut--shaped regions given in Table \ref{regions1}. We consider the scaled regions $c_j D = D_j$ where we now take a decreasing scaling parameter $c_j = (1-j/10)$. Notice that for the three regions 
$$|c_j D| = \frac{3\pi}{4} c_j^2 \,\, \text{for the kite--shaped,} \quad |c_j D| =\frac{809\pi}{3200} c_j^2 \,\, \text{for the star--shaped}$$
and
$$|c_j D| =\frac{5\pi}{8} c_j^2 \,\, \text{for the peanut--shaped regions}.$$ 
As in the previous convergence study, we approximate the rate of convergence via \eqref{RoC-formula}. From the scaling Theorem \ref{eig-scaling} we expect that the rate of convergence for these regions should be given by $p=1$. 
\begin{table}[H]
\centering
\begin{tabular}{r|r|r|r|r|r|r}
$c_j$ & $k_1($Kite$)$ & rate $p_j$ & $k_1($Star$)$ & rate $p_j$ & $k_1($Peanut$)$ & rate $p_j$\\
\hline
1.0 & $1.9166 $ & ---           & 3.2672  & ---       & 2.1309 & ---\\
0.9 & $2.1296$  & 1.0000 & 3.6302  & 1.0000  & 2.3677 & 1.0000\\
0.8 & $ 2.3958$ & 1.0000 & 4.0839  & 1.0000  & 2.6637 & 1.0000\\
0.7 & $2.7381$ & 1.0000 & 4.6674   & 1.0000  & 3.0442 & 1.0000\\
0.6 & $3.1944$ & 1.0000 & 5.4453   & 1.0000  & 3.5515 & 1.0000\\
0.5 & $3.8333$ & 1.0000 & 6.5343   & 1.0000  & 4.2619 & 1.0000\\
0.4 & $4.7916$ & 1.0000 & 8.1679   & 1.0000  & 5.3273 & 1.0000\\
0.3 & $6.3889$ & 1.0000 & 10.8905  & 1.0000 & 7.1031 & 1.0000\\
\hline
\end{tabular}
 \caption{\label{table4}The first CTE and the rate of convergence for the scaled regions where $c_j = (1-j/10)$ for $j=0, \cdots 7$ for the kite--shaped, star--shaped and peanut--shaped regions.}
\end{table}
The result for this convergence study is given in Table \ref{table4}. As the theory dictates, we see that for the scaled regions we have that the rate of convergence $p=1$. Just as for the CTE for a disk, we see that the lower bound in Theorem \ref{eig-smallvolumelimit} is sharp for the case of a scaled region. This is due to the fact that the scaled regions are shape regular as the measure decreases.

Another numerical test we give here is to see how the first CTE depends on the perimeter of the region. From \eqref{firsteig-disk} we see that the CTEs are decreasing for the case when $D$ is a disk. But this is deceptive since the perimeter decreases if and only if the area decreases for the case of a disk. We will now consider regions that have the same area but different perimeter. To this end, we define the the region 
\begin{align} 
\partial D_\varepsilon = \big(0.75\cos(t)+\varepsilon\cos(2t),\sin(t) \big)^\top  \quad \text{for $0<\varepsilon \ll 1$ with $t \in [0,2\pi)$} \label{regions3}
\end{align}
where we have that $|D_\varepsilon | = \frac{3\pi}{4}$ for any appropriate $\varepsilon$ but the perimeter $\ell(\varepsilon)$ is given by 
$$\ell(\varepsilon) = \int_0^{2\pi}\sqrt{(0.75\sin(t)+2\varepsilon \sin(2t))^2+\cos^2(t)}\,\mathrm{d}t.$$
In Table \ref{table5}, we report the first CTE for various values of $\varepsilon$. Here, we see that the eigenvalues seem to be increasing with respect to the perimeter. It is open whether there is a monotonicity relation with respect to the perimeter. Here, we compute $\ell(\varepsilon)$ via standard numerical integration. 
\begin{table}[!ht]
\centering
 \begin{tabular}{r|c|c|c|c|c}
                                     $\varepsilon$ & $0.1$ & $0.2$ & $0.3$ & $0.4$ & $0.5$ \\
\hline
$k_1(D_\varepsilon)$  & 1.88515 & 1.89716 & 1.91665 & 1.94279 & 1.97439\\
$\ell(\varepsilon)$       & 5.56836 & 5.69900 & 5.92478 & 6.24952 & 6.66757\\
\hline
\end{tabular}
\caption{\label{table5}The first real CTE for the region $D_\varepsilon$ given in \eqref{regions3} using $\varepsilon\in \{0.1,0.2,0.3,0.4,0.5\}$ and the corresponding perimeter $\ell(\varepsilon)$.}
\end{table}

With this, we have provided an extensive numerical study of the CTEs. We have shown that the theory developed in Sections \ref{RQ-analysis} and \ref{eigpair-analysis} is valid. From our numerics, we have also considered the questions: (1) what is an appropriate upper bound for $k^2_1(D)$ with respect to $|D|$ and (2) are the CTEs monotone with respect to the perimeter i.e. length of $\partial D$ provided the measure is fixed. These are interesting questions that can be studied as part of future work on the clamped transmission eigenvalues.

\section{Conclusion}\label{fin}
In this paper, we studied the clamped transmission eigenvalue problem associated with the flexural wave scattering of a plane wave by a clamped obstacle. The direct scattering problem has applications to non--destructive testing of thin elastic materials. Here we established fundamental analytical results of the principal eigenpair. By introducing an equivalent variational formulation \eqref{varform-4RQ} and a suitable Rayleigh quotient \eqref{eqn_frac_form}, we have rigorously proven several key results for this eigenvalue problem. The results that are proven include: 
\begin{itemize}
\item The principal eigenvalue is monotone decreasing with respect to the region.
\item The corresponding eigenfunction is of a single sign.
\item A Faber--Krahn type inequality holds for the first clamped transmission eigenvalue.
\item The principal eigenvalue is strictly bounded by the first Dirichlet eigenvalue.
\end{itemize}

These theoretical results were substantiated by numerical investigations, which also explored the asymptotic behavior of the eigenvalues as the domain size tends to zero. Future work could focus on establishing a sharp upper bound for the principal eigenvalue and investigating its monotonicity with respect to the perimeter for a fixed measure. Another interesting study could be the approximation of the clamped transmission eigenvalue via a FEM. In particular, in \cite{scatterpoles-imped,MaSun,resonances-wDtN} the authors give a rigorous study of the approximation of scattering poles via a FEM as well as boundary integrals. Since the scattering poles are also eigenvalues for a differential operator on an unbounded domain, the methods studied in the aforementioned manuscript could be applicable to the clamped transmission eigenvalue problem. Also, the method of fundamental solutions was used in \cite{MFS} to compute the interior transmission eigenvalues associated with acoustic scattering. One can also consider the inside--outside duality (see,  e.g., \cite{armin,inout-iso}) to further connect this eigenvalues to the associated direct scattering problem.  \\

\noindent{\bf Acknowledgments:} The research of author I. Harris was partially supported by the NSF DMS Grant 2509722. \\


\end{document}